%%%%%%%%%%%%%%%%%%%%%%%%%
%%  FICHIER PLAIN TEX  %%
%%%%%%%%%%%%%%%%%%%%%%%%%

                  %%%%%%%%%%%%%%%%%%%%
                  %%  COMMENTAIRES  %%
                  %%%%%%%%%%%%%%%%%%%%
%
% Pour les th\'eor\`emes, propositions, lemmes,
% corollaires, assertions, conjectures, probl\`emes,
% ... et d\'efinitions, en Fran\c{c}ais ou en Anglais,
% les taper en "petites capitales" (small capitals)
% et leurs textes en italique, sauf pour les
% d\'efinitions o\`u il est en romain, avec,
% \'eventuellement, l'objet ou le concept \`a
% d\'efinir en italique. Les faire pr\'ec\'eder
% par un "\medskip" et suivre par un "\smallskip".
% Les remarques, exemples, cas etc... sont \`a
% taper en italique et leurs textes en romain.
%
% Exemple :
%
% \medskip
%
% {\sc Th\'eor\`eme~2.1.~--~}{\it \'Enonc\'e du
% th\'eor\`eme.}
%
% \smallskip
%
% {\it Remarque\/}~2.2.~--~}Texte de la remarque
%                     

%%%%%%%%%%%%%%%%%%%%%%%%%%%%%%%%%%%%%%%%%%%%%%%%%%%%%%%%%%%%%%%
%%%%%%%%%%%%%%%%%%%%%%%%%%%%%%%%%%%%%%%%%%%%%%%%%%%%%%%%%%%%%%%

\magnification=1000
\overfullrule=0pt

\font\eightrm=cmr8
\font\sixrm=cmr10 at 6pt
\textfont0=\tenrm
\scriptfont0=\sevenrm
\scriptscriptfont0=\fiverm
\def\rm{\fam0\tenrm}

%%%%%%%%%%%%%%%%
%%  FONT CMM  %%
%%%%%%%%%%%%%%%%

% \font\sixi=cmmi6
% \font\eighti=cmmi8
% \font\ninei=cmmi9
% \font\twelvei=cmmi12
% \font\fourteeni=cmmi10 at 14pt
\font\teni=cmmi10
\font\seveni=cmmi7
\font\fivei=cmmi5
\textfont1=\teni
\scriptfont1=\seveni
\scriptscriptfont1=\fivei
\def\mit{\fam1}
\def\oldstyle{\fam1\teni}

%%%%%%%%%%%%%%%%%
%%  FONT CMSY  %%
%%%%%%%%%%%%%%%%%

\font\fivesy=cmsy6 at 5pt
\font\sixsy=cmsy6
\font\sevensy=cmsy7
\font\eightsy=cmsy8

\font\tensy=cmsy10

%%%%%%%%%%%%%%%%%
%%  FONT CMTI  %%
%%%%%%%%%%%%%%%%%

\font\eightit=cmti8

%%%%%%%%%%%%%%%%%%%
%%  FONT CMBFTI  %%
%%%%%%%%%%%%%%%%%%%

\font\bfti=cmbxti10

%%%%%%%%%%%%%%%%%%%%%%%%%%%%%%%%%%%%%
%%  FONT CMCSC (Petites capitales) %%

%%%%%%%%%%%%%%%%%%%%%%%%%%%%%%%%%%%%%
\font\sc=cmcsc10

%%%%%%%%%%%%%%%%%
%%  FONT MSBM  %%
%%%%%%%%%%%%%%%%%

\font\sixbboard=msbm7 at 6pt
\font\eightbboard=msbm8 % msbm10 at 8pt
\font\eightbboard=msbm8 % msbm10 at 8pt
\font\tenbboard=msbm10
\font\sevenbboard=msbm7
\font\fivebboard=msbm5
\newfam\bboardfam
\textfont\bboardfam=\tenbboard
\scriptfont\bboardfam=\sevenbboard
\scriptscriptfont\bboardfam=\fivebboard
\def\bb{\fam\bboardfam\tenbboard}
\let\oldbb=\bb
\def\bb #1{{\oldbb #1}}

%%%%%%%%%%%%%%%%%%%%%%%%%%%%%
%%  FONT EUFM  (Gothique)  %%
%%%%%%%%%%%%%%%%%%%%%%%%%%%%%

\font\tengoth=eufm10
\font\sevengoth=eufm7
\font\fivegoth=eufm5
\newfam\gothfam
\textfont\gothfam=\tengoth
\scriptfont\gothfam=\sevengoth
\scriptscriptfont\gothfam=\fivegoth
\def\goth{\fam\gothfam\tengoth}

%%%%%%%%%%%%%%%%%
%%  FONT CMBX  %%
%%%%%%%%%%%%%%%%%

\font\fourteenbf=cmbx10 at 14pt

\font\tenbf=cmbx10

\font\eightbf=cmbx8
\font\sevenbf=cmbx7
\font\sixbf=cmbx5 at 6pt
\font\fivebf=cmbx5
\font\tenbf=cmbx10
\font\sevenbf=cmbx7
\font\fivebf=cmbx5 
\newfam\bffam
\textfont\bffam=\tenbf
\scriptfont\bffam=\sevenbf
\scriptscriptfont\bffam=\fivebf
\def\bf{\fam\bffam\tenbf}

%MSSYMB.TeX     
\catcode`\@=11
\font\tenmsx=cmr8 %msxm10
\font\sevenmsx=cmr8 %msxm7
\font\fivemsx=cmr8 %msxm5
\font\tenmsy=cmr8 %msym10
\font\sevenmsy=cmr8 %msym7
\font\fivemsy=cmr8 %msym5
\newfam\msxfam
\newfam\msyfam
\textfont\msxfam=\tenmsx  \scriptfont\msxfam=\sevenmsx
\scriptscriptfont\msxfam=\fivemsx
\textfont\msyfam=\tenmsy  \scriptfont\msyfam=\sevenmsy
\scriptscriptfont\msyfam=\fivemsy

%%%%%%%%%%%%%%%%%%%%%%%%%%%%%%%%%%%%%%%%%%%%%%%%%%%%%%%%%%%%

%%%%%%%%%%%%%%%%%%%%%%%%%%%%%%%%%%%%%%%%%%%%%%%%%%%%%%%%%%%%
\font\eightrm=cmr8
\font\eighti=cmmi8
\font\eightsy=cmsy8
\font\eightbf=cmbx8

%\font\eighttt=cmtt8
\font\eightit=cmti8

%\font\eightsl=cmsl8
\font\sixrm=cmr10 at 6pt % cmr6
\font\sixi=cmmi6
\font\sixsy=cmsy6
\font\sixbf=cmbx5 at 6pt %\font\sixbf=cmbx6

%\skewchar\eighti='177
\skewchar\sixi='177

%\skewchar\eightsy='60 \skewchar\sixsy='60

%%%%%%%%%%%%%%%%%%%%%%%%%%%%%%%%%%%%%%%%%%%%%%%%%%%%%%%%%%%%
\font\sl=cmsl10
\font\eightsl=cmsl8
\font\tt=cmr8 %cmstt10
\font\eighttt=cmtt8
\def\tenpoint{%
\textfont0=\tenrm \scriptfont0=\sevenrm \scriptscriptfont0=\fiverm
\def\rm{\fam0\tenrm}%
\textfont1=\teni \scriptfont1=\seveni \scriptscriptfont1=\fivei
\def\mit{\fam\@ne}\def\oldstyle{\fam1\teni}%
\textfont2=\tensy \scriptfont2=\sevensy \scriptscriptfont2=\fivesy
\def\itfam{4}\textfont\itfam=\tenit
\def\it{\fam\itfam\tenit}%
\def\slfam{5}\textfont\slfam=\tensl
\def\sl{\fam\slfam\tensl}%
\def\bffam{6}\textfont\bffam=\tenbf \scriptfont\bffam=\sevenbf
\scriptscriptfont\bffam=\fivebf
\def\bf{\fam\bffam\tenbf}%
\def\ttfam{7}\textfont\ttfam=\tentt
\def\tt{\fam\ttfam\tentt}%
\textfont\bboardfam=\tenbboard \scriptfont\bboardfam=\sevenbboard
\scriptscriptfont\bboardfam=\fivebboard

\def\bb{\fam\bboardfam\tenbboard}%

%  \font\sc=cmcsc10
\abovedisplayskip=6pt plus 2pt minus 6pt
\abovedisplayshortskip=0pt plus 3pt
\belowdisplayskip=6pt plus 2pt minus 6pt
\belowdisplayshortskip=7pt plus 3pt minus 4pt
\smallskipamount=3pt plus 1pt minus 1pt
\medskipamount=6pt plus 2pt minus 2pt
\bigskipamount=12pt plus 4pt minus 4pt
\normalbaselineskip=12pt
\setbox\strutbox=\hbox{\vrule height8.5pt depth3.5pt width0pt}%
\normalbaselines\rm}

\def\eightpoint{%

\textfont0=\eightrm \scriptfont0=\sixrm \scriptscriptfont0=\fiverm
\def\rm{\fam0\eightrm}%
\textfont1=\eighti \scriptfont1=\sixi \scriptscriptfont1=\fivei
\def\oldstyle{\fam1\eighti}%
\textfont2=\eightsy \scriptfont2=\sixsy \scriptscriptfont2=\fivesy
\textfont\slfam=\eightit
\def\sl{\fam\itfam\eightit}%
\textfont\slfam=\eightsl
\def\sl{\fam\slfam\eightsl}%
\textfont\bffam=\eightbf \scriptfont\bffam=\sixbf
\scriptscriptfont\bffam=\fivebf
\def\bf{\fam\bffam\eightbf}%
\textfont\ttfam=\eighttt
\def\tt{\fam\ttfam\eighttt}%
\textfont\bboardfam=\eightbboard \scriptfont\bboardfam=\sixbboard
\scriptscriptfont\bboardfam=\fivebboard
\def\bb{\fam\bboardfam\eightbboard}%

%  \font\sc=cmcsc8

\abovedisplayskip=9pt plus 2pt minus 6pt
\abovedisplayshortskip=0pt plus 2pt
\belowdisplayskip=9pt plus 2pt minus 6pt  
\belowdisplayshortskip=5pt plus 2pt minus 3pt  
\smallskipamount=2pt plus 1pt minus 1pt
\medskipamount=4pt plus 2pt minus 1pt
\bigskipamount=9pt plus 3pt minus 3pt
\normalbaselineskip=9pt
\setbox\strutbox=\hbox{\vrule height7pt depth2pt width0pt}%
\normalbaselines\rm}

%\input macro
%%%%%%%%%%%%%%%%%%%%%%%%%%%%%%%%%%%%%%%%%%%%%%%%%%%%%%%%%%%%
%%  HAUT-DE-PAGE (EX. : AUTEUR COURANT ET TITRE COURANT)  %%
%%  BAS-DE-BAGE                                           %%
%%%%%%%%%%%%%%%%%%%%%%%%%%%%%%%%%%%%%%%%%%%%%%%%%%%%%%%%%%%%

\newif\ifpagetitre           \pagetitretrue
\newtoks\hautpagetitre       \hautpagetitre={\hfil}
\newtoks\baspagetitre        \baspagetitre={\hfil}

\newtoks\auteurcourant       \auteurcourant={\hfil}
\newtoks\titrecourant        \titrecourant={\hfil}

\newtoks\hautpagegauche
         \hautpagegauche={\hfil\the\auteurcourant\hfil}
\newtoks\hautpagedroite
         \hautpagedroite={\hfil\the\titrecourant\hfil}

\newtoks\baspagegauche
         \baspagegauche={\hfil\tenrm\folio\hfil}
\newtoks\baspagedroite
         \baspagedroite={\hfil\tenrm\folio\hfil}

\headline={\ifpagetitre\the\hautpagetitre
\else\ifodd\pageno\the\hautpagedroite
\else\the\hautpagegauche\fi\fi }

\footline={\ifpagetitre\the\baspagetitre
\global\pagetitrefalse
\else\ifodd\pageno\the\baspagedroite
\else\the\baspagegauche\fi\fi }

%%%%%%%%%%%%%%%%%%%%%%%%%%%%%%%%%%%%%%%%%%%%%%%%%%%%%%%%%%%%%%%

           %%%%%%%%%%%%%%%%%%%%%%%%%%%%%%%%%%%%%%%%%%%
           %%  UN SEMBLANT DES GUILLEMETS FRANCAIS  %%
           %%%%%%%%%%%%%%%%%%%%%%%%%%%%%%%%%%%%%%%%%%%

\def\og{\leavevmode\raise.3ex
     \hbox{$\scriptscriptstyle\langle\!\langle$~}}
\def\fg{\leavevmode\raise.3ex
     \hbox{~$\!\scriptscriptstyle\,\rangle\!\rangle$}}

\def\ogg{\leavevmode\raise.3ex
     \hbox{$\scriptstyle\scriptscriptstyle\langle\!\langle$~}}
\def\fgg{\leavevmode\raise.3ex
     \hbox{~$\!\scriptstyle\scriptscriptstyle\,\rangle\!\rangle$}}

%%%%%%%%%%%%%%%%%%%%%%%%%%%%%%%%%%%%%%%%%%%%%%%%%%%%%%%%%%%%%%%%%%%

      %%%%%%%%%%%%%%%%%%%%%%%%%%%%%%%%%%%%%%%%%%%%%%%%%%%%%%
      %%  TRUC AVEC MACHIN AU-DESSUS ET BIDULE AU-DESSOUS %%
      %%                (MODE MATHEMATIQUE)               %%
      %%%%%%%%%%%%%%%%%%%%%%%%%%%%%%%%%%%%%%%%%%%%%%%%%%%%%%

\def\build #1_#2^#3{\mathrel{\mathop{\kern 0pt #1}
     \limits_{#2}^{#3}}}

%%%%%%%%%%%%%%%%%%%%%%%%%%%%%%%%%%%%%%%%%%%%%%%%%%%%%%%%%%%%%%%%%%%

%%%%%%%%%%%%%%%%%%%%%%%%%%%%%%%%%%%%%%%%%%%%%%%%%%%%%%%%%%%%%%%%%%%%

         %%%%%%%%%%%%%%%%%%%%%%%%%%%%%%%%%%%%%%%%%%%%%%%%%%%
         %% FLECHES HORIZONTALES ET VERTICALES AVEC TRUC  %%
         %%   AU-DESSUS ET MACHIN AU-DESSOUS (MODE MATH)  %%
         %%%%%%%%%%%%%%%%%%%%%%%%%%%%%%%%%%%%%%%%%%%%%%%%%%%

\def\hfl[#1][#2][#3]#4#5{\kern -#1
 \sumdim=#2 \advance\sumdim by #1 \advance\sumdim by #3
  \smash{\mathop{\hbox to \sumdim {\rightarrowfill}}
   \limits^{\scriptstyle#4}_{\scriptstyle#5}}
    \kern-#3}

\def\vfl[#1][#2][#3]#4#5%
 {\sumdim=#1 \advance\sumdim by #2 \advance\sumdim by #3
  \setbox1=\hbox{$\left\downarrow\vbox to .5\sumdim {}\right.$}
   \setbox1=\hbox{\llap{$\scriptstyle #4$}\box1
    \rlap{$\scriptstyle #5$}}
     \vcenter{\kern -#1\box1\kern -#3}}

%%%%%%%%%%%%%%%%%%%%%%%%%%%%%%%%%%%%%%%%%%%%%%%%%%%%%%%%%%%%%%%%%%%%

                   %%%%%%%%%%%%%%%%%%%%%%%%%%%%%
                   %% DIAGRAMMES COMMUTATIFS  %%
                   %%%%%%%%%%%%%%%%%%%%%%%%%%%%%
\newdimen\sumdim
\def\diagram#1\enddiagram
    {\vcenter{\offinterlineskip
      \def\tvi{\vrule height 10pt depth 10pt width 0pt}
       \halign{&\tvi\kern 5pt\hfil$\displaystyle##$\hfil\kern 5pt
        \crcr #1\crcr}}}

%% SYNTAXE

%%  $$
%%  \diagram
%%   ... & & ... & \cr
%%   ... & & ... & \cr
%%   .................
%%   ... & & ... & \cr
%%  \enddiagram
%%  $$        

%%%%%%%%%%%%%%%%%%%%%%%%%%%%%%%%%%%%%%%%%%%%%%%%%%%%%%%%%%%%%%%%%%%%%%%%%%%%%%%%%%%%%%
%%%%%%%%%%%%%%%%%%%%%%%%%%%%%%%%%%%%%%%%%%%%%%%%%%%%%%%%%%%%%%%%%%%%%%%%%%%%%%%%%%%%%%

%\def\Doteq{\build{=}_{\scriptscriptstyle\bullet}^{\scriptscriptstyle\bullet}}

%%%%%%%%%%%%%%%%%%%%%%%%%%%%%%%%%%%%%%%%%%%%%%%%%%%%%%%%%%%%%%%%%%%%

               %%%%%%%%%%%%%%%%%%%%%%%%%%%%%%%%%%%
               %% ESPACE VERTICAL EN MODE MATH  %%
               %%%%%%%%%%%%%%%%%%%%%%%%%%%%%%%%%%%

\def\vspace[#1]{\noalign{\vskip #1}}

%%%%%%%%%%%%%%%%%%%%%%%%%%%%%%%%%%%%%%%%%%%%%%%%%%%%%%%%%%%%%%%%%%%%

                      %%%%%%%%%%%%%%%%%%%
                      %%  ENCADREMENT  %%
                      %%%%%%%%%%%%%%%%%%%

\def\boxit [#1]#2{\vbox{\hrule\hbox{\vrule
     \vbox spread #1{\vss\hbox spread #1{\hss #2\hss}\vss}%
      \vrule}\hrule}}

  %% INDIQUE LA FIN D'UNE PREUVE
 %% OPERATEUR DE LA CHALEUR
 %% OPERATEUR DE LA CHALEUR
                           %% EN INDICE OU EN EXPOSANT

%%%%%%%%%%%%%%%%%%%%%%%%%%%%%%%%%%%%%%%%%%%%%%%%%%%%%%%%%%%%%%%%%%%%

\catcode`\@=11
\def\Eqalign #1{\null\,\vcenter{\openup\jot\m@th\ialign{
    \strut\hfil$\displaystyle{##}$&$\displaystyle{{}##}$\hfil
    &&\quad\strut\hfil$\displaystyle{##}$&$\displaystyle{{}##}$
    \hfil\crcr #1\crcr}}\,}
\catcode`\@=12

\catcode`\@=11
\def\displaylinesno #1{\displ@y\halign{
 \hbox to\displaywidth{$\@lign\hfil\displaystyle ##\hfil$}&
  \llap{\rm ##} \crcr #1\crcr}}
  
\def\ldisplaylinesno #1{\displ@y\halign{
 \hbox to\displaywidth{$\@lign\hfil\displaystyle ##\hfil$}&
  \kern-\displaywidth\rlap{\rm ##}%
   \tabskip\displaywidth\crcr #1\crcr}} 
\catcode`\@=12

\vsize= 20.5cm
\hsize=15cm
\hoffset=3mm
\voffset=3mm

%%%%%%%%%%%%%%%%%%%%%%%%%%%%%%%%%%%%%%%%%%%%%%%%
%% DEBUT DE DEFINITIONS PROPRES A CE FICHIER  %%
%%%%%%%%%%%%%%%%%%%%%%%%%%%%%%%%%%%%%%%%%%%%%%%%

\def \frak {\goth}
\def \mathfrak {\goth}

\def \mathcal {\cal}
\def \ad {\mathop {\rm ad}}
\def \adj {\mathop {\rm adj}}
\def \tr {\mathop {\rm tr}}

\def \mapright#1{\smash{\mathop{\longrightarrow}\limits^{#1}}}

%%%%%%%%%%%%%%%%%%%%%%%%%%%%%%%%%%%%%%%%%%%%%%%
%% FIN DE DEFINITIONS PROPRES A CE FICHIER   %%
%%%%%%%%%%%%%%%%%%%%%%%%%%%%%%%%%%%%%%%%%%%%%%%

%%%%%%%%%%%%%%%%%%%%%%%%%%%%%%%%%%%%%%%%%%%%%%%%%%%%%%%%%%%%%%%%%%%
%%%%%%%%%%%%%%%%%%%%%%%%%%%%%%%%%%%%%%%%%%%%%%%%%%%%%%%%%%%%%%%%%%%

%%%%%%%%%%%%%%%%%%%%%%%%%%%%%%%%%%%%%%%%%%%%%
%%  AUTEUR(S) COURANT(S) SUR PAGES PAIRES  %%
%%%%%%%%%%%%%%%%%%%%%%%%%%%%%%%%%%%%%%%%%%%%%

\auteurcourant{\eightbf M.~Andler, A.~Dvorsky, S.~Sahi
\hfill}

%%%%%%%%%%%%%%%%%%%%%%%%%%%%%%%%%%%%%%%%
%%  TITRE COURANT SUR PAGES IMPAIRES  %%
%%%%%%%%%%%%%%%%%%%%%%%%%%%%%%%%%%%%%%%%

\titrecourant{\hfill\eightbf Deformation quantization and
invariant distributions}

%%%%%%%%%%%%%%%%%%%%%%%%%%%%%%%%%%%%%%%%%%%%%%%%%%%%%%%%%%%%%%%

{\parindent=0pt\eightbf
%C. R. Acad. Sci. Paris, t.~, S\'erie~I,
%p.~000--000, 1999

%%%%%%%%%%%%%%%%%
%% RUBRIQUE(S) %%
%%%%%%%%%%%%%%%%%

%G\'eom\'etrie Diff\'erentielle/\hskip -.5mm{\eightbfti Differential
%Geometry}

%(Th\'eorie des groupes/\hskip -.5mm{\eightbfti Group Theory})

\vskip 1.5cm

%%%%%%%%%%%%%%%%%%%%%%%%%%%%%%%%%%%%%%%%%%%%%%%%%%%%%%%%%%%%%%%%%%%

%%%%%%%%%%%%%%%%%%%%%%%%%%%%%%%%%%%%%%%%
%%  TITRE IN EXTENSO EN ANGLAIS  %%
%%%%%%%%%%%%%%%%%%%%%%%%%%%%%%%%%%%%%%%%

{\fourteenbf
Deformation quantization and invariant distributions
}
%%%%%%%%%%%%%%%%%%%%%%%%%%%%%%%%%%%%%%%%%%%
%%  MENTION "NOTE PR\'ESENT\'EE PAR ...  %%
%%%%%%%%%%%%%%%%%%%%%%%%%%%%%%%%%%%%%%%%%%%
\parindent=-1mm\footnote{}{
%Note pr\'esent\'ee par Isra\"el G{\sixbf ELFAND.}
}

\vskip 15pt

%%%%%%%%%%%%%%%%%%%%%%%%%%%%%%%%%%%
%%  Pr\'enoms NOMS (DES AUTEURS) %%
%%%%%%%%%%%%%%%%%%%%%%%%%%%%%%%%%%%

\noindent{\bf
Martin ANDLER~${}^{\hbox{\sevenbf a}}$,
Alexander DVORSKY~${}^{\hbox{\sevenbf b}}$,
Siddhartha SAHI~${}^{\hbox{\sevenbf c}}$
}
\vskip 6pt

%%%%%%%%%%%%%%%%%%%%%%%%%%%%%%%%%%%%%%%%%%%%%%%%%%
%%  ADRESSES DES AUTEURS ET COURRIEL OU E-MAIL  %%
%%%%%%%%%%%%%%%%%%%%%%%%%%%%%%%%%%%%%%%%%%%%%%%%%%

{\parindent=3mm
\item{${}^{\hbox{\fivebf a}}$}D\'{e}partement de Mathematiques,
Universit\'{e} de
Versailles Saint Quentin, \ 78035 Versailles C\'{e}dex.
Courrier \'electronique~: andler@math.uvsq.fr \par}

\smallskip

{\parindent=3mm
\item{${}^{\hbox{\fivebf b, c}}$}Mathematics Department, Rutgers
University, New Brunswick, NJ 08903.\hfill\par
e-mail~: dvorsky@math.rutgers.edu, sahi@math.rutgers.edu
\par}
\medskip

%%%%%%%%%%%%%%%%%%%%%%%%%%%%%%%%%%%%%%%%%%%%%%%%%%%%%
%%  MENTION "(RE\C{C}U LE ..., ACCEPT\'E LE ...)"  %%
%%%%%%%%%%%%%%%%%%%%%%%%%%%%%%%%%%%%%%%%%%%%%%%%%%%%%

%(Re\c{c}u le , accept\'e apr\`es r\'evision le )
}

\vskip3pt

\line{\hbox to 2cm{}\hrulefill\hbox to 1.5cm{}}

\vskip3pt

%%%%%%%%%%%%%%%%
%%  ABSTRACT  %%
%%%%%%%%%%%%%%%%

{\parindent=2cm\rightskip 1.5cm\baselineskip=9pt
\item{\bf Abstract.~}{\eightrm
We study Kontsevich's deformation quantization for the dual of a
finite dimensional Lie algebra $\frak g$. Regarding elements of
$\cal S(\frak g)$ as distributions on $\frak g$, we show that
the $\star$-multiplication operator ($r\mapsto r\star p$) is
a differential operator with analytic germ at $0$. We use this
to establish a conjecture of Kashiwara and Vergne which, in turn,
gives a new proof of Duflo's result on the local solvability of
bi-invariant differential operators on a Lie group. ~}
%{\eightrm\copyright~Acad\'emie des Sciences/Elsevier, Paris}
\par}

\vskip 10pt

%%%%%%%%%%%%%%%%%%%%%%%%
%%  TITRE EN FRANCAIS  %%
%%%%%%%%%%%%%%%%%%%%%%%%

{\parindent=2cm\rightskip 1.5cm
{\bfti Quantification par d\'eformation et distributions invariantes}\par}

\vskip 10pt

%%%%%%%%%%%%%%%%%%
%%  R\'ESUM\'E  %%
%%%%%%%%%%%%%%%%%%

{\parindent=2cm\rightskip 1.5cm\baselineskip=9pt
\item{\bf R\'esum\'e.~}{\eightrm
Nous \'etudions la quantification par d\'eformations dans le cas du dual
d'une alg\`ebre de Lie de dimension finie $\frak g$.
Consid\'erant les \'el\'ements de $\cal S(\frak g)$ comme des distributions
de support ${0}$ sur $\frak g$, nous montrons que
la $\star$-multiplication \`a droite ($r \mapsto r \star p$)
est un op\'erateur diff\'erentiel
\`a germes analytiques en ${0}$. Nous utilisons ceci pour \'etablir une
conjecture de Kashiwara-Vergne, ce qui fournit une nouvelle
d\'emonstration du th\'eor\`eme de Duflo sur la
r\'esolubilit\'e locale des op\'erateurs diff\'erentiels bi-invariants sur
un groupe de Lie.
%~\copyright~Acad\'emie des Sciences/Elsevier, Paris
}
\par}

\vskip 10pt
\line{\hbox to 2cm{}\hrulefill\hbox to 1.5cm{}}
\medskip

\parindent=3mm
\vskip 5mm

%%%%%%%%%%%%%%%%%%%%%%%%%%%%%%%%
%%   VERSION FRANCAISE ABREGEE %%
%%%%%%%%%%%%%%%%%%%%%%%%%%%%%%%%
{\bf Version fran\c caise abr\'eg\'ee}

\medskip

Soit $G$ un groupe de Lie r\'eel connexe de dimension finie
d'alg\`ebre de Lie $\frak{g}$, et soit
$\mathcal{U}=\mathcal{U(\frak{g})}$  et $\mathcal{S=S(\frak{g})}$,
respectivement
l'alg\`ebre enveloppante et l'alg\`ebre sym\'etrique de $\frak{g}$, que
nous identifions
\`a l'alg\`ebre de convolution des distributions de support $1$ sur le
groupe
$G$, et \`a
l'alg\`ebre de convolution des distributions de support $0$ sur $\frak g$.
Notons
$\mathcal{I}=\mathcal{S}^{\mathcal{\frak{g}}},
\mathcal{Z}=\mathcal{U}^{\mathcal{\frak{g}}}$
les sous-alg\`ebres d'invariants correspondantes.

Soit ${\bf U}$ et ${\bf  S}$ les espaces de {\it germes} en $1$
et $0$ de distributions sur $G$ et $\frak{g}$, et ${\bf  Z}$ et ${\bf  I}$
les sous-espaces d'invariants
correspondants. Nous avons les diagrammes commutatifs

$$ 
\matrix{
\cal U  &\mapright{}    &{\bf U}         \cr
\uparrow        &       &\uparrow\cr
\cal Z  &\mapright{}    &{\bf  Z}        \cr
}
\qquad \qquad
\matrix{
\cal S  &\mapright{}    &{\bf  S}        \cr
\uparrow        &       &\uparrow\cr
\cal I  &\mapright{}    &{\bf  I}        \cr
}$$

Soit $\exp_{\ast} : {\bf  S} \to {\bf U}$ et
$\exp^{\ast} : {\bf U} \to {\bf  S}$ l'image directe et
inverse au sens des (germes de)
distributions. Soit $q$ la fonction analytique sur un voisinage
de $0$ dans $\frak{g}$~:
$$q(x):=\sqrt{\det\left( e^{\ad(x)/2}-e^{-\ad(x)/2}\over \ad(x)\right)}.$$
Soit $\eta$ l'application de $S$ dans $U$
$$\eta(p)=\exp_{\ast}(p\ q)$$

Utilisant la m\'ethode des orbites de Kirillov, Duflo [Du] a montr\'e que la
restriction de $\eta$ \`a
$\cal I$ est un isomorphisme d'{\it alg\`ebres} de $\cal I$ sur $\cal Z$,
c'est-\`a-dire~:
$$\eta(p_1 \ast_{\frak g} p_2) = \eta(p_1) \ast_G \eta(p_2) \quad {\rm pour
}\; p_1,p_2 \in \cal{I},$$
o\`u $\ast_G$ et $\ast_{\frak g}$ sont respectivement le produit de
convolution sur le groupe $G$
et l'alg\`ebre de Lie $\frak{g}$.

Soit $\frak{D}$ l'alg\`ebre des germes en $0$ d'op\'erateurs
diff\'erentiels \`a
coefficients analytiques sur $\frak{g}$. Soit $\frak R$ l'id\'eal \`a droite
de $\frak D$ engendr\'e par les germes de champs de vecteurs
$\adj_a$ d\'efinis par l'action adjointe de $G$ sur $\frak g$ :
$$\adj{}_{a}(f)(x)={d \over dt}f(\exp(-ta)\cdot x)|_{t=0},\; a\in
\frak{g}.$$
Naturellement, $\bf S$ est un $\frak D$-module \`a droite --
nous noterons donc l'action de $\frak D$ sur les distribution
\`a droite.

Nous d\'emontrons le r\'esultat suivant, conjectur\'e en
g\'en\'eral, et d\'emontr\'e dans le cas r\'esoluble, par Kashiwara --Vergne
[KV, Corollary 1].
\medskip
{\sc Theor\`eme~1. -- }{\it Pour $p$
dans $\mathcal{I}$, l'op\'erateur $D_{p}: {\bf  S} \to {\bf  S}$ appartient
\`a $\frak R$~:
$$P\cdot D_{p} = \exp^{\ast}\left(\eta(P\ast_{\frak g}
p)-\eta(P)\ast_G\eta(p)\right).$$
}
\smallskip
Comme cela est expliqu\'e dans [KV], ce r\'esultat a plusieurs
cons\'equence
importantes. D'abord, comme les \'el\'ements de
$\frak{R}$ annulent les distributions dans
${\bf  I}$ nous obtenons~:
\medskip
{\sc Theor\`eme~2. -- }{\it Pour tout $p$ dans $\mathcal{I}$ et $P$ dans
${\bf  I}$, nous avons
$$\eta(P\ast_{\frak g} p)=\eta(P)\ast_G\eta(p).$$}
\smallskip
Rappelons ensuite qu'une distribution $Q$ sur $G$ est appel\'ee
distribution propre s'il existe un caract\`ere
$\chi:\mathcal{Z}\longrightarrow\bf C$ tel que
$$Q\ast_G z=\chi(z)Q \ \ {\rm pour\  tout}\ z\ {\rm dans}\ \mathcal{Z}.$$
On d\'efinit de m\^eme les distributions propres
sur l'alg\`ebre de Lie; nous obtenons le
\medskip

{\sc Theor\`eme~3. -- }{\it L'application $\eta$ envoie les germes de
distributions propres et invariantes sur l'alg\`ebre de Lie sur ceux
du groupe.} \smallskip
Ces r\'esultats ont des cons\'equences importantes pour l'analyse des
op\'erateurs diff\'erentiels bi-invariants sur $G$ (voir version en anglais
ci-dessous).

Nous d\'emontrons le Th\'eor\`eme 1 en nous appuyant sur le
travail remarquable de M. Kontsevich sur la formalit\'e du complexe de
Hochschild pour
une  vari\'et\'e diff\'erentiable [Ko]. Les \'etapes de la d\'emonstration
sont les
suivantes. On d\'emontre d'abord que, dans notre cas, on peut poser $h=1$
dans le produit
$\star$, puis que, pour $p$ fix\'e dans $\cal S(\frak g)$, l'application
$\partial^\star_p:r\in {\cal S} \mapsto r\star p$
appartient \`a $\frak D$. Nous d\'emontrons ensuite qu'une certaine s\'erie
formelle
consid\'er\'ee par Kontsevich [Ko] converge dans un voisinage de $0$ dans
$\frak g$
vers une  fonction $\tau(x)$. On d\'emontre ensuite que pour $p \in \cal I$
l'op\'erateur $T_p = \partial_p \tau - \tau \partial^\star_{p\tau}$
appartient \`a $\frak R$, o\`u $\partial_p$ est l'op\'erateur
de convolution par $p$ dans $\frak g$. On compare enfin $T_p$ et $D_p$
en d\'emontrant la formule $D_p = T_p \tau^{-1} q.$
\bigskip

%%%%%%%%%%%%%%%%%%%%%%%%%%%%%%%%%%%%%%%%%%%%%%%%%%%%%%%%%%%%%
%%%%%%%%%%%%%%%%%%%%%%%%%%%%%%%%%%%%%%%%%%%%%%%%%%%%%%%%%%%%%

\vskip 6pt
\centerline{\hbox to 2cm{\hrulefill}}
\vskip 9pt

%%%%%%%%%%%%%%%%%%%%%%%%%%%%%%%%%%%%%%%%%%%%%%%%%%%%%%%%%%%%%
%%%%%%%%%%%%%%%%%%%%%%%%%%%%%%%%%%%%%%%%%%%%%%%%%%%%%%%%%%%%%

%%%%%%%%%%%%%%%%%%%%%%%
%%  TEXTE PRINCIPAL EN ANGLAIS  %%
%%%%%%%%%%%%%%%%%%%%%%%

\noindent{\bf 1. Invariant distributions}
\medskip

Let $G$ be a finite dimensional real Lie group with Lie algebra
$\frak{g}$. We write $\mathcal{U}=\mathcal{U(\frak{g})},
\mathcal{S=S(\frak{g})}$ for the enveloping and symmetric
algebras of $\frak{g}$, and denote by
$\mathcal{I}=\mathcal{S}^{\mathcal{\frak{g}}},
\mathcal{Z}=\mathcal{U}^{\mathcal{\frak{g}}}$ the
subalgebras of $\frak{g}$-invariants.

Let ${\bf U}$ and ${\bf S}$ be the spaces of {\it germs} at $1$ and
$0$ of distributions on $G$ and $\mathcal{\frak{g}}$, and
let ${\bf  Z}$ and ${\bf  I}$ be the subspaces of $\frak{g}$-invariant
distributions.
We identify $\mathcal{U}$ and $\mathcal S$ with the convolution
algebras of distributions with point support $1\in G$ and
$0\in\frak{g}$. Then we have the following inclusions

$$ 
\matrix{
\cal U  &\mapright{}    & {\bf U}        \cr
\uparrow        &       &\uparrow\cr
\cal Z  &\mapright{}    &{\bf  Z}        \cr
}
\qquad \qquad
\matrix{
\cal S  &\mapright{}    &{\bf S}         \cr
\uparrow        &       &\uparrow\cr
\cal I  &\mapright{}    &{\bf I}         \cr
}$$

We denote by $\exp_\ast: {\bf S}\longrightarrow {\bf U}$ and $\exp^\ast:
{\bf U}\longrightarrow {\bf S}$ the pushforward and pullback of
distributions (germs)
under the exponential map; and define $\eta:{\bf S}\longrightarrow {\bf U}$
by
$$\eta(p)=\exp_\ast(pq)$$
where $q$ is the following analytic function on $\frak g$~:
$$q(x):=\sqrt{\det\left(  e^{\ad(x)/ 2}
-e^{-\ad(x)/ 2}\over \ad(x)\right)}.$$

Using the orbit method initiated by Kirillov, Duflo [Du] proved that the
restriction of $\eta$ to $\mathcal{I}$ is an {\it algebra} isomorphism
from $\mathcal{I}$ to $\mathcal{Z}$, i.e.
$$ \eta(p_1\ast_{\frak g}p_2)=\eta(p_1)\ast_G\eta(p_2)
\quad {\rm for }\; p_1,p_2 \in \cal{I}$$
where $\ast_{\frak g}$ and $\ast_G$ denote
convolution in the Lie algebra and Lie group respectively.

Let $\frak{D}$ be the algebra of germs at $0$ of analytic coefficient
differential operators on $\frak{g}$, and let
$\frak{R}$ be the right ideal of $\frak{D}$ generated by
the germs of the adjoint vector fields
$$\adj{}_{a}(f)(x)={d \over dt}f(\exp(-ta)\cdot x)|_{t=0},\; a\in
\frak{g}.$$

By virtue of the pairing between functions and distributions, ${\bf S}$ is
naturally
a {\it right\/} $\frak{D}$-module, and to emphasise this fact we shall write
the action on distributions on the right.

We shall establish the following result which was conjectured by Kashiwara
and Vergne, and proved by them for solvable groups [KV, Corollary 1].

\medskip
{\sc Theorem~1. -- }{\it For $p$ in $\mathcal{I}$,
the following differential operator $D_{p}:{\bf S}\rightarrow {\bf S}$
lies in $\frak{R}$:}
$$P\cdot D_{p}=\exp^{\ast}\big(\eta(P\ast_{\frak g}
p)-\eta(P)\ast_G\eta(p)\big).$$
\smallskip
As is explained in [KV], this result has several important consequences.
First, since operators from $\frak{R}$ annihilate distributions from ${\bf
I}$,
we obtain the following extension of the Duflo isomorphism~:
\medskip
{\sc Theorem~2. -- }{\it For $p$ in $\mathcal{I}$ \ and $P$ in ${\bf I}$, we
have
$$
\eta(P\ast_{\frak g} p)=\eta(P)\ast_G\eta(p).$$}
\smallskip

Next, recall that a distribution $Q$ on $G$ is called an
eigendistribution if there is a
character $\chi:\mathcal{Z\longrightarrow}\bf C$ such that
$$Q\ast_G z=\chi(z)Q \ {\rm for}\;{\rm all}\ z\ {\rm in}\ \mathcal{Z}{\rm
}$$
Eigendistributions on the Lie algebra are defined similarly, and we
obtain the following:
\medskip
{\sc Theorem~3. -- }{\it The map $\eta$ takes germs of invariant
eigendistributions on the Lie algebra to those on the Lie group.}
\smallskip

Another consequence is the following important result,
which was conjectured by L. Schwartz and proved by Duflo in [Du]
by different methods. (See also [He], [Ra], [DR] for various special cases):
\medskip
{\sc Corollary. -- }{\it Every nonzero bi-invariant differential operator
on $G$ is locally solvable.}
\smallskip

This follows immediately from Theorem 2 and Rais' lemma from [Ra] on
the existence of {\it invariant} fundamental solutions on the
Lie algebra $\frak{g}$.
\bigskip
\bigskip
\noindent{\bf 2. Deformation quantization}
\medskip
We will deduce Theorem 1 from careful analysis of M. Kontsevich's
remarkable result on the formality of the Hochschild complex of a smooth
manifold [Ko].

Kontsevich's main result leads to the construction of an explicit $\star
$-product on the algebra of functions on an arbitrary smooth Poisson
manifold $(X,\gamma)$. This $\star$-product is given by a formal power
series of bidifferential operators (depending on a parameter $\hbar$).

In the special situation when $X= \frak{g}^{\ast}$ is the dual of a finite
dimensional Lie algebra $\frak{g}$, and $\gamma$ is the usual Poisson
bracket (dual to the Lie bracket); we can set $\hbar=1$ and Kontsevich's
$\star$-product descends to an actual product on the algebra of polynomial
functions on $\frak{g}^{\ast}$, which can be naturally identified with
$\mathcal{S=S(\frak{g})}$.

If, as before, we consider $\mathcal{S}$ as the algebra of distributions
with point support (at $0\in\mathcal{\frak{g}}$), then the $\star$-product
can be viewed as a new convolution on $\mathcal{S}$. The advantage of the
distributional point of view is the following key result:
\medskip
{\sc Theorem 4. -- }{\it Given $p\in\mathcal{S}$ there is a unique
element $\partial_{p}^{\star}$ in
$\frak{D}$ such that for all $r\in\mathcal{S}$
$$
r\star p=r\cdot \partial_{p}^{\star}.$$}
\smallskip
(As explained in the previous section, we are writing differential operators
on the right.)

The proof of this result involves examining the Kontsevich $\star$-product
formula, replacing multiplications by differentiation and vice versa. This
results in an expression for $\partial_{p}^{\star}$ as an infinite sum of
polynomial coefficient differential operators of {\it bounded\/} order.
The coefficients of the sum are given by certain integrals which, while
hard to compute, can be readily estimated to deduce Theorem 4. The linearity
of
$\gamma$ is crucial in these considerations.

Next, in [Ko, 8.3.3] Kontsevich introduces a certain formal power series
$$ S_{1}(x)=\exp\left(  \sum_{k=1}^{\infty}c_{2k}\tr(\ad x)^{2k}\right),$$
where the constants $c_{2k}$ (denoted $c_{2k}^{(1)}$ in [Ko]) are
expressed as integrals of certain smooth differential forms (``angle
forms'') over certain compact manifolds with corners (compactified
configuration
spaces, introduced in [FM]). It can be checked that
$$ c_{n}\sim\alpha^{n}$$
for some positive constant $\alpha$. Hence in a neighborhood of $0$ the
power series $S_{1}(x)$ converges to an analytic function which we shall
denote
by $\tau(x).$
\medskip
{\sc Theorem 5. -- }{\it For $p \in \cal I$, the operator
$T_{p}=\partial_{p}\tau-\tau\partial_{p\tau }^{\star}$ belongs to
$\frak{R}$.}
\smallskip
Here $\partial_{p\tau}^{\star}$ is the differential operator defined in
Theorem 4 and $r\cdot\partial_{p}:=r\ast_{\frak g} p$.

To establish this result, we study the expression
$$(r\ast_{\frak g} p)\tau-(r\tau)\star (p\tau)$$
for $r\in \cal S$ and $p \in \cal I$, as in section 8.2 of [Ko].
This enables us to express $T_p$ as an infinite sum of differential
operators each belonging to $\frak{R}$, and Theorem 5
follows from Theorem 4.

Theorem 1 is now a consequence of the following lemma:
\medskip
{\sc Lemma 1. -- }{\it
For $p\in \cal S$, we have
$$ D_{p}= T_{p}\tau^{-1}q.$$
}
\smallskip

To prove this lemma it is enough to show that
$r\cdot D_p=r\cdot(T_{p}\tau^{-1}q)$ for $r$ in $\cal S$.
In turn, this follows from the fact
that there is a canonical isomorphism
$\kappa:\mathcal{S}\rightarrow \cal U$
(denoted $I_{alg}^{-1}$ in [Ko]) which
satisfies the relations [Ko, 8.3.1]
$$\kappa(r\star p)=\kappa(r)\ast_G\kappa(p)$$
and [Ko, 8.3.4]
$$\kappa(r\tau)=\eta(r)=\exp_\ast(rq)\ {\rm for}\ r\in\mathcal{S}.$$

\bigskip

{\it Acknowledegment. --} The authors would like to thank M. Kontsevich for
sharing his insights with us, especially for explaining the proof of Theorem
8.2 of [Ko]. Moreover, this paper would not have been written without the
help of I. Gelfand and F. Knop, and that of the other participants of the
Gelfand seminar at Rutgers, especially R. Wilson, C. Weibel, J. Barros-Neto,
M. Nathanson, J. Lagarias, D.\thinspace Stone, O. Stoyanov, and A.
Kazarnovski-Krol,
who helped us tremendously in the enjoyable task of discussing, dissecting,
and
digesting Kontsevich's formidable manuscript.
\bigskip

%%%%%%%%%%%%%%%%%%%%%%%
%%   BIBLIOGRAPHIE   %%
%%%%%%%%%%%%%%%%%%%%%%%

\centerline{\bf Bibliographic References}

\medskip

{\baselineskip=9pt\eightrm

\item{[Du]} M. Duflo, Op\'{e}rateurs diff\'{e}rentiels
bi-invariants sur un groupe de Lie, Ann. Sc. \'{E}c. Norm. Sup 10
(1977), 267-288
\smallskip

\item{[DR]} M. Duflo, M. Rais, Sur l'analyse harmonique
sur les groupes de Lie r\'{e}solubles, Ann. Sc. \'{E}c. Norm. Sup 9 (1976),
107-114.
\smallskip

\item{[FM]} W. Fulton, R. MacPherson, Compactification of
configuration spaces, Ann. of Math. 139 (1994), 183-225.
\smallskip

\item {[He]} S. Helgason, The surjectivity of invariant
differential operators on a symmetric space, Ann. of Math. 98
(1973), 451-479.

\item{[KV]} M. Kashiwara, M. Vergne, The Campbell-Hausdorff
formula and invariant hyperfunctions, Inventiones Math. 47 (1978), 249-272.
\smallskip

\item{[Ko]} M. Kontsevich, Deformation quantization of
Poisson manifolds I, e-print (math.QA/9709040).
\smallskip

\item{[Ra]} M. Rais, Solutions \'{e}l\'{e}mentaires des
op\'{e}rateurs diff\'{e}rentiels bi-invariants sur un groupe de Lie
nilpotent, C. R. Acad. Sci. 273 (1971), 495-498.
\smallskip
}
\end